\documentclass[11pt,a4paper]{amsart}
\date{\today}   %%
\usepackage{hyperref}
\usepackage{latexsym,amsmath,amsfonts,amscd,amssymb,graphics}
\usepackage[all]{xy}
\usepackage[utf8]{inputenc}
\usepackage{tikz-cd}
\usepackage[margin=1in]{geometry}

\theoremstyle{plain}  % default
\newtheorem{theorem}{Theorem}[section]
\newtheorem*{theorem*}{Theorem}

\newtheorem{lemma}[theorem]{Lemma}
\newtheorem{proposition}[theorem]{Proposition}

\theoremstyle{definition}

\newtheorem*{notation*}{Notation}
\newtheorem{definition}[theorem]{Definition}

\theoremstyle{remark}
\newtheorem{example}[theorem]{Example}
\newtheorem{remark}[theorem]{Remark}

\newtheorem*{claim*}{Claim}

\numberwithin{equation}{section}

\renewcommand{\leq}{\leqslant}

\renewcommand{\geq}{\geqslant}

\renewcommand{\setminus}{\smallsetminus}

\newcommand{\R}{\mathbb{R}}

\newcommand{\Z}{\mathbb{Z}}

\newcommand{\Pin}{\mathrm{Pin}}

\newcommand{\PGL}{\mathrm{PGL}}

\newcommand{\PSO}{\mathrm{PSO}}

\newcommand{\GL}{\mathrm{GL}}

\newcommand{\Or}{\mathrm{O}}
\newcommand{\PO}{\mathrm{PO}}
\newcommand{\Map}{\mathrm{Map}}

\newcommand{\CW}{\mathrm{CW}}

\DeclareMathOperator{\Hom}{Hom}

\newcommand{\Aut}{\operatorname{Aut}}

\newcommand{\smtrx}[1]{\left (\begin{smallmatrix}#1\end{smallmatrix}\right)}

\let\oldmarginpar\marginpar
\renewcommand\marginpar[1]{\oldmarginpar{\tiny\bf\begin{flushleft} #1
\end{flushleft}}}

\begin{document}

%%%%%%%%%%%%%%%%%%%%%%%%%%%%%%%%%%%%%%%%%%%%%%%%%%%%%%%%%%%%%%%%%%%
\title[Bundles on $2$-dimensional $\CW$-complexes]{Principal bundles on $2$-dimensional $\CW$-complexes with disconnected structure group}
%%%%%%%%%%%%%%%%%%%%%%%%%%%%%%%%%%%%%%%%%%%%%%%%%%%%%%%%%%%%%%%%%%%

\author[A. Oliveira]{Andr\'e Oliveira}
\address{}
\email{}

\thanks{
%%%%%%%%%%%%%%%%%%%%%%%%%%%%%%%%%%%%%%%%%%%%%%%%%%%%%%%%%%%%%%%%%%
The author is grateful to Gustavo Granja (Instituto Superior Técnico, Lisboa) for his generous support and for introducing him to these subjects, to Peter Gothen for useful discussions, and also to the anonymous referee for the comments provided.
%%%%%%%%%%%%%%%%%%%%%%%%%%%%%%%%%%%%%%%%%%%%%%%%%%%%%%%%%%%%%%%%%%
Author partially supported by CMUP (UIDB/00144/2020) with national funds.
}


\begin{abstract}
Given any topological group $G$, the topological classification of principal $G$-bundles over a finite $\CW$-complex $X$ is long-known to be given by the set of free homotopy classes of maps from $X$ to the corresponding classifying space $BG$. This classical result has been long-used to provide such classification in terms of explicit characteristic classes. However, even when $X$ has dimension $2$, there is a case in which such explicit classification has not been explicitly considered. This is the case where $G$ is a Lie group, whose group of components acts non-trivially on its fundamental group $\pi_1G$.  
Here we deal with this case and obtain the classification, in terms of characteristic classes, of principal $G$-bundles over a finite $\CW$-complex of dimension $2$, with $G$ is a Lie group such that $\pi_0G$ is abelian.
\end{abstract}

\maketitle

%\tableofcontents

\section{Introduction}

It is a classical fact that, for a path-connected $\CW$-complex $X$ and for a topological group $G$, the set of topological types of principal $G$-bundles over $X$ is in bijection with the set $[X,BG]$ of free homotopy classes of maps from $X$ to the classifying space $BG$ of $G$. It is, however, many times useful to have a more explicit classification of such bundles, that is, to have a more detailed description of the set $[X,BG]$, for example in terms of characteristic classes.

In this note we provide such a description in a setting which is quite classical, and hence which ought to be very well-known. We obtain an explicit and complete classification of $G$-bundles over a $2$-dimensional connected finite $\CW$-complex $X$, for any Lie group $G$ with $\pi_0G$ abelian\footnote{The abelian condition for $\pi_0G$, besides being naturally verified for most Lie groups, also makes several of the technical arguments much easier to handle with. For example, under that condition, the classifying space $B\pi_0G$ is a topological abelian group (cf. Sections \ref{sec:4.1} and \ref{sec:5.2}).}, in terms of cohomology classes. This includes two classical cases. Firstly, if $G$ is connected, such bundles are classified by the cohomology group $H^2(X,\pi_1G)$. Secondly, there is a natural action of $\pi_0G$ on $\pi_1G$, and if this action is trivial, then $G$-bundles over $X$ are classified by the product $H^1(X,\pi_0G)\times H^2(X,\pi_1G)$. 
Nonetheless, such $\pi_0G$-action can be non-trivial and, somehow surprisingly, this case does not seem to have been explicitly treated so far. 
An instance where this phenomenon occurs is the case of the projective orthogonal group of even rank, which we explicitly deal with as an example. 

The importance of the topological classification of principal $G$-bundles over any finite, $2$-dimensional, $\CW$-complex by itself or by its connections with  other areas of Mathematics and Physics is obvious, so the mentioned untreated case deserves to be considered. This is the purpose of this article. Just to give an example where such a topological classification plays an important role, consider the following. Suppose $G$ is semisimple and consider the space $\mathcal{R}(X,G)=\mathrm{Hom}^{red}(\pi_1X,G)/G$ of reductive representations $\rho:\pi_1X\to G$ modulo the $G$-action by conjugation. These spaces, or spaces deeply related to them, arise from several different natural points of view (Physics, Gauge Theory, Hyperbolic Geometry, etc.) and have been intensively studied in the past decades, mostly in the case $X$ is a closed oriented surface, where deep connections arise with the theory of holomorphic vector bundles \cite{narasimhan-seshadri:1965} and also with the theory of Higgs bundles, under the so-called non-abelian Hodge correspondence (there are many references for this subject: the seminal paper is \cite{hitchin:1987} and an overview and other references may be found for instance in \cite{wentworth:2016}). However, the spaces $\mathcal{R}(X,G)$ have also been studied for compact, non-orientable surfaces \cite{ho-liu:2003,baird:2010-I,baird:2010-II,schaffhauser:2011} and even for any finite $2$-dimensional $\CW$-complex \cite{daskalopoulos-mese-wilkin:2018}. Now, the topological type of $G$-bundles over $X$ provides a way to distinguish certain connected components of  $\mathcal{R}(X,G)$ and to give a lower bound on the number of such components (for example, it is known that if $G$ is a complex Lie group and $X$ is an orientable closed surface, such bound is exact \cite{garcia-prada-oliveira:2017}). Hence, the theorem of this article gives a way to achieve this for any Lie group $G$ with $\pi_0G$ abelian.

To give a hint on the classification theorem, note first that if a $G$-bundle over $X$ is given by the homotopy class of a map $X\to BG$, then the homotopy class of composite $X\to BG\to B\pi_0G$ defines a topological invariant $\mu_1\in [X,B\pi_0G]\cong H^1(X,\pi_0G)$ of the given $G$-bundle, measuring the obstruction to reduce its structure group to the identity component of $G$. Having this, the classification result (see Theorem \ref{thm:class-thm}) states that there is a bijection between the subset of $[X,BG]$ consisting of isomorphism classes of principal $G$-bundles over $X$ with invariant $\mu_1\in H^1(X,\pi_0G)$ and the quotient set $H^2(X,\pi_1\mathcal{G}_{\mu_1})/\pi_0G$. Here $\pi_1\mathcal{G}_{\mu_1}$ denotes the local system obtained from $\mu_1$ and from a natural action of $\pi_0G$ in $\pi_1G$ and where $\pi_0G$ acts on $H^2(X,\pi_1\mathcal{G}_{\mu_1})$ also via the same action of $\pi_0G$ in $\pi_1G$.

\section{The classification and examples}
\subsection{A first topological invariant}
Let $X$ be a path-connected $2$-dimensional $\CW$-complex. Let also $G$ be a Lie group, with (discrete) group of connected components $\pi_0G$.
Recall that isomorphism classes of principal $G$-bundles over $X$ are topologically classified by the set $[X,BG]$ of free homotopy classes of maps from $X$ to the classifying space $BG$. There is a natural map
\begin{equation}\label{maptoBpi0}
\chi:[X,BG]\to[X,B\pi_0G], \ \ \ \chi([f])=[p_{0,*}\circ f],
\end{equation}
where $p_{0,*}:BG\to B\pi_0G$ is the map induced from the canonical projection $p_0:G\to\pi_0G$.

From here one readily defines a first topological invariant of a $G$-bundle $E$ over $X$.
\begin{definition}\label{mu1}
Let $E$ be a $G$-bundle over $X$ represented by a map $f:X\to BG$. Define $$\mu_1(E)=\chi([f])\in [X,B\pi_0G].$$ 
\end{definition}

Since $B\pi_0G$ is the Eilenberg-Maclane space $K(\pi_0G,1)$, it follows that $[X,B\pi_0G]\cong H^1(X,\pi_0G)$. $\mu_1(E)$ is obviously a topological invariant of $E$ and it represents the isomorphism class of the $\pi_0G$-bundle on $X$ obtained from $E$ through the projection $G\to\pi_0G$. 

From now on we fix an element $\mu_1\in [X,B\pi_0G]$ and our task is to classify $G$-bundles $E$ on $X$ such that $\mu_1(E)=\mu_1$. 

\subsection{The classification theorem and examples}
Notice that the action $$\Psi:G\longrightarrow\Aut(G)$$ of $G$ on itself by conjugation induces an action of $\pi_0G$ on the higher homotopy groups $\pi_iG$ (by $\pi_iG$ we always mean homotopy groups based at the identity of $G$): for $a\in G$, consider the induced automorphism $\Psi(a)_*:\pi_iG\to\pi_iG$ and, if $b\in G$ is in the same component as $a$, then $\Psi(a)$ is homotopic to $\Psi(b)$ via a path joining $a$ and $b$, hence $\Psi(a)_*=\Psi(b)_*$.

In particular, we are interested in the action of $\pi_0G$ on $\pi_1G$, which we will denote by 
\begin{equation}\label{action pi0G in pi1G}
 \Psi(-)_*:\pi_0G\longrightarrow\Aut(\pi_1G).
\end{equation}

Assume now, once and for all, that the group $\pi_0G$ is abelian and consider the homomorphism in $\pi_1$, $\mu_{1*}:\pi_1X\to\pi_0G$, induced from the invariant $\mu_1:X\to B\pi_0G$ that we have previously fixed (recall that $\pi_0G=\pi_1(B\pi_0G)$). Given this, we have an induced action of $\pi_1X$ on $\pi_1G$ 
\begin{equation}\label{Psimu1}
\Psi(-)_*\circ \mu_{1*}:\pi_1X\longrightarrow\Aut(\pi_1G)
\end{equation}
and so this action defines a local system $\pi_1\mathcal{G}_{\mu_1}$ on $X$.

Let $H^*(X,\pi_1\mathcal{G}_{\mu_1})$ be the cohomology of $X$ with values in the local system $\pi_1\mathcal{G}_{\mu_1}$. Recall that this is the cohomology of the cochain complex $C^*_{\Z[\pi_1X]}(\widetilde{X},\pi_1G)=\Hom_{\Z[\pi_1X]}(C_*(\widetilde{X}),\pi_1G)$, where $\Z\pi=\Z[\pi_1X]$ is the group ring of $\pi_1X$ and $C_*(\widetilde{X})$ is the $\Z\pi$-module of chains on the universal cover $\widetilde X$ of $X$.

Note now that  $\Psi(-)_*$ in \eqref{action pi0G in pi1G} induces, by composition, a $\pi_0G$-action on the $\Z\pi$-module $C^2_{\Z[\pi_1X]}(\widetilde{X},\pi_1G)$. 
Explicitly, $a\in \pi_0G$ acts on a $2$-cochain $\tau\in C^2_{\Z[\pi_1X]}(\widetilde{X},\pi_1G)$ as 
\begin{equation}\label{action pi0G on C2}
a\cdot\tau=\Psi(a)_*\circ\tau. 
\end{equation}
Since $\widetilde{X}$ has dimension $2$, every $2$-cochain is a cocycle hence, passing to the quotient, \eqref{action pi0G on C2} yields a $\pi_0G$-action on $H^2(X,\pi_1\mathcal{G}_{\mu_1})$.
Let $$H^2(X,\pi_1\mathcal{G}_{\mu_1})/\pi_0G$$ be the corresponding quotient set. 

We are now ready to state the result concerning the topological classification of principal $G$-bundles over $X$.

\begin{theorem}\label{thm:class-thm}
 Let $G$ be a Lie group with $\pi_0G$ abelian and let $X$ be a $2$-dimensional connected $\CW$-complex. 
There is a bijection between the set of isomorphism classes of continuous principal $G$-bundles over $X$ with invariant $\mu_1\in H^1(X,\pi_0G)$ and the (non-empty) quotient set $H^2(X,\pi_1\mathcal{G}_{\mu_1})/\pi_0G$.
\end{theorem}

This means that the map $$\chi:[X,BG]\to H^1(X,\pi_0G)$$ defined in \eqref{maptoBpi0} is such that there is bijection of sets $\chi^{-1}(\mu_1)\simeq H^2(X,\pi_1\mathcal{G}_{\mu_1})/\pi_0G$.

\begin{remark}
Of course if one has a topological space $M$ which is homotopically equivalent to a finite $2$-dimensional $\CW$-complex $X$, then $[M,BG]\cong[X,BG]$, so the above theorem also applies for the classification of principal $G$-bundles over $M$. For instance, $M$ could be the complement of a finite number of points in a compact $3$-dimensional manifold. 
\end{remark}

\begin{remark}
Using different methods, Theorem \ref{thm:class-thm} has actually been proved before in \cite{oliveira:2011}, in the particular case $X$ is a closed oriented surface.
\end{remark}

We will prove the theorem in section \ref{sec: proof}. In the next sections, we will briefly recall some notions (of obstruction theory and of Postnikov sections) which will be required for it. Before that, let us now consider  some explicit examples (most of them quite classical) of applications of this theorem.

\begin{example} Let $X$ be any finite connected $2$-dimensional $\CW$-complex and $G$ a Lie group with $\pi_0G$ abelian.
 \begin{enumerate}
 \item If $G$ is connected, then $G$-bundles over $X$ are classified by $$[X,BG]\cong H^2(X,\pi_1G)$$ as it is well-known (for surfaces, see for example \cite{ramanathan:1975}, Proposition $5.1$).
 \item If $\pi_0G$ acts trivially on $\pi_1G$, then $G$-bundles over $X$ are classified by $$[X,BG]\cong H^1(X,\pi_0G)\times H^2(X,\pi_1G).$$
Because $\pi_0G$ is abelian, this case yields
$$[X,BG]\cong\pi_0G^{2g}\times\pi_1G$$
when $X$ is a closed (i.e., compact and without boundary) oriented surface of genus $g$, and yields
$$[X,BG]\cong(\pi_0G)_2\times\pi_0G^{k-1}\times\pi_1G/2\pi_1G$$
when $X$ is a closed non-orientable surface which is a connected sum of $k$ copies of the projective plane $\R\mathbb{P}^2$. Here $(\pi_0G)_2$ denotes the $2$-torsion subgroup of $\pi_0G$.
\item\label{ex:O(n)} When $G=\Or(n)$, the group of orthogonal transformations of $\R^n$, with $n\geq 3$, then we are in the situation of the previous item and the classification is given by the first and second Stiefel-Whitney classes $w_1,w_2$.

\item If $X$ is simply-connected, then principal $G$-bundles over $X$ are classified by $$[X,BG]\cong H^2(X,\pi_1G)/\pi_0G.$$
In particular, this holds for the $2$-sphere $S^2$ and gives $[S^2,BG]\cong \pi_1G/\pi_0G$. Again this is quite classical: see \cite{steenrod:1999}, Section $18$.

\item\label{ex:O(2)} Here is another classical example, in which $\pi_0G$ acts non-trivially on $\pi_1G$, and so which we can deduce from the above theorem. Consider $G=\Or(2)$, the group of orthogonal transformations of $\R^2$.

Then $\pi_0\Or(2)=\Z_2$ acts non-trivially on $\pi_1\Or(2)=\Z$, by changing the sign of the generator. Indeed, a generator can be represented by the loop $\gamma(\theta)=\smtrx{\cos(2\pi\theta) &\sin(2\pi\theta) \\ -\sin(2\pi\theta) & \cos(2\pi\theta)}$, with $\theta\in[0,1]$, and the non-trivial element of $\pi_0\Or(2)$, which we can represent by  $\smtrx{1 & 0 \\ 0 & -1}$, takes it to the inverse loop $\gamma(-\theta)$. 

\begin{enumerate}
\item
Suppose $X$ is a closed oriented surface of genus $g$. Then we have the isomorphism $H^2(X,\Z)\cong\Z$ by cap product with the fundamental class. Since the action of $\pi_0\Or(2)=\Z_2$ on $H^2(X,\Z)$ by \eqref{action pi0G on C2} is by post-composition on the coefficients, the isomorphism $H^2(X,\Z)\cong\Z$ is $\pi_0\Or(2)$-equivariant with the action on the right-hand-side given by \eqref{action pi0G in pi1G}. Note that this argument always applies for any coefficients $\pi_1G$, and hence applies generally and not just for $G=\Or(2)$.

Hence, $\Or(2)$-bundles over $X$ with $\mu_1=0$, i.e. $w_1=0$ (hence which reduce to $\mathrm{SO}(2)\cong S^1$-bundles), are classified by $H^2(X,\Z)/\Z_2\cong\Z/\Z_2\cong\Z_{\geq 0}$. 
%To see what this is in a more explicit way, note that such bundles are of the form $L\oplus L^{-1}$ for $L$ a line bundle on $X$. One knows that actually the first Chern class $c_1(L)\in H^2(X,\Z)$ of $L$ distinguishes them. Hence actually, $H^2(X,\Z)/\Z_2\cong H^2(X,\Z)$.
%In case $X$ is a closed oriented surface, this gives the degree of $L$.

Consider now $\Or(2)$-bundles over $X$ with first Stiefel-Whitney class $w_1=\mu_1\neq 0$ (so non-orientable rank $2$ orthogonal bundles on $X$). These are hence classified by $H^2(X,\Z_{w_1})/\Z_2$.
Then, by Poincaré duality with local coefficients (see \cite{spanier:1993}, Theorem 10.4) and by \cite[Proposition 5.14 (1)]{davis-kirk:2001} (or \cite[Theorem 3.2, VI]{whitehead:1978}),
\begin{equation}\label{eq:H2 O2}
H^2(X,\Z_{w_1})\cong H_0(X,\Z_{w_1})\cong\Z/H_{w_1}=\pi_1\Or(2)/H_{w_1},
\end{equation} where 
\begin{equation}\label{eq:Hw1}
H_{w_1}=\langle x-\gamma\cdot x\,|\, x\in\pi_1\Or(2), \gamma\in\pi_1X\rangle
\end{equation} and $\gamma\cdot x$ is the action of $\gamma$ on $x$ determined by $w_1$ and by \eqref{Psimu1}. From the calculation at the beginning of the example, $H_{w_1}\cong 2\Z$, so $H^2(X,\Z_{w_1})\cong\Z_2$. Now we have to look at the action of $\Z_2$ on the set $H^2(X,\Z_{w_1})$ induced by \eqref{action pi0G on C2}. As above, the isomorphism $H^2(X,\Z_{w_1})\cong\Z_2=\{0,1\}$ is equivariant with respect to the $\pi_0\Or(2)$-actions \eqref{action pi0G on C2} and  \eqref{action pi0G in pi1G}, because duality is also given by cap product and also because the isomorphism of \cite{davis-kirk:2001} Proposition 5.14 (1) is $\pi_0G$-equivariant. Again this holds in general and not just for $G=\Or(2)$.
But the $\Z_2=\pi_0\Or(2)$-action on $\Z_2=\{0,1\}$ changes the sign of the generator, hence preserves the parity, thus is trivial. So $H^2(X,\Z_{w_1})/\Z_2\cong\Z_2$.

We conclude then the following common knowledge fact: rank $2$ orthogonal bundles over a closed oriented surface $X$ of genus $g$ are topologically classified by the characteristic classes in
$$(\mu_1,\mu_2)\in(\{0\}\times\Z_{\geq 0})\cup ((\Z_2)^{2g}\setminus\{0\})\times\Z_2$$
where $\mu_1$ is the first Stiefel-Whitney class and $\mu_2$ is the (non-negative) degree when $\mu_1=0$ and the second Stiefel-Whitney class when $\mu_1\neq 0$.
\item
Suppose now that $X$ a closed but non-orientable surface. Then $H^2(X,\Z)\cong \Z/2\Z\cong\Z_2$, and we have to see how $\Z_2=\pi_0\Or(2)$ acts via \eqref{action pi0G on C2} on $H^2(X,\Z)$. This is the same as considering the action \eqref{action pi0G in pi1G} of $\Z_2=\pi_0\Or(2)$ on $\Z/2\Z\cong\Z_2$ because the isomorphism $H^2(X,\Z)\cong \Z/2\Z$ is equivariant by arguments similar to the orientable case (just take instead the fundamental class of $X$ in homology with local coefficients, cf. \cite{hatcher:2002}, Example 3H.3); this holds in general and not only for $\Or(2)$. From above, such action changes the sign of the generator, hence preserves the parity. So $H^2(X,\Z)/\Z_2\cong(\Z/2\Z)/\Z_2\cong\Z/2\Z\cong\Z_2$.

Suppose now that $\mu_1=w_1\neq 0$. Write 
\begin{equation}\label{eq:wX}
w_X:\pi_1X\to\Aut(\Z)\cong\Z_2
\end{equation} for the corresponding non-trivial orientation character associated to the oriented double cover (this is really the first Stiefel-Whitney class of $X$ i.e. of its tangent bundle), and $\Z_{w_X}$ for the associated local system on $X$. 
Then, as above, using Poincaré duality  (see \cite[Theorem 10.4]{spanier:1993}), and \cite[Proposition 5.14 (1)]{davis-kirk:2001}, we find that
\begin{equation}\label{eq:H2nonorienmu1neq0}
H^2(X,\Z_{w_1})\cong H_0(X,\Z_{w_X}\otimes_{\Z[\pi_1X]}\Z_{w_1})\cong(\Z/H_{w_X}\otimes_{\Z}\Z/H_{w_1})\cong \Z_2\otimes_{\Z}\Z_2\cong\Z_2,
\end{equation}
where $H_{w_1}$ is given analogously to \eqref{eq:Hw1} and $H_{w_X}=\langle x-\gamma\cdot x\,|\, x\in\Z, \gamma\in\pi_1X\rangle$ with $\pi_1X$ acting by the non-trivial orientation character $w_X$ in \eqref{eq:wX} (thus so that $H_{w_X}\cong 2\Z$). 
Here we used again that the two first isomorphisms are $\pi_0G$-equivariant with respect to \eqref{action pi0G on C2} and  \eqref{action pi0G in pi1G} (and this is true for any group, not just $\Or(2)$).

Concluding, rank $2$ orthogonal bundles over a non-orientable surface $X$ which is a connected sum of $k$ copies of $\R\mathbb{P}^2$, are topologically classified by
$$(\mu_1,\mu_2)\in(\{0\}\times\Z_2)\cup (\Z_2^k\setminus\{0\}\times\Z_2)=\Z_2^k\times\Z_2.$$ 
\end{enumerate}

 \item\label{ex:PO(n)} Let us see now a perhaps slightly less known case. We will use, without reference, some of the results we used in the previous example, which hold in general, such as the various versions of Poincaré duality and the $\pi_0G$-equivariance of the isomorphisms from the various $H^2$ under the actions \eqref{action pi0G in pi1G} and \eqref{action pi0G on C2}.

 Let $n\geq 4$ be even and consider the projective orthogonal group $\PO(n)=\Or(n)/\Z_2$, where $\Z_2=\{\pm I_n\}$ is the center of $\Or(n)$. (If $n$ is odd $\PO(n)$ is connected.) Since $n$ is even, $\pi_0\PO(n)\cong\Z_2$ and 
$$\pi_1\PO(n)=\begin{cases}
    \Z_2\times\Z_2 & \text{if } n=0\ \text{mod}\ 4 \\
    \Z_4 & \text{if } n=2\ \text{mod}\ 4.
  \end{cases}$$
More precisely, the universal cover of $\PO(n)$ is the 
$\Pin(n)$ and, if $p:\Pin(n)\to\PO(n)$ is the covering
  projection, then, as a set (and using the abelian notation) $\ker(p)=\{0,1,\omega_n,-\omega_n\}$ where $\omega_n=e_1\cdots e_n$ is the oriented volume element of $\Pin(n)$ in the standard construction of this group via the Clifford algebra $\mathrm{Cl}(n)$ (see, for example, \cite{lawson-michelson:1989}). So, if $n=0 \ \text{mod}\ 4$, $\pm\omega_n$ are elements of order $2$, while if $n=2 \ \text{mod}\ 4$, $\pm\omega_n$ have order $4$. 
 It turns out that $\pi_0\PO(n)$ acts on $\pi_1\PO(n)$ by leaving $0$ and $1$ fixed and identifying $\pm\omega_n$. More precisely, recall that $\Pin(n)$ is a group with two connected components, $\Pin(n)^-$ and $\mathrm{Spin}(n)$, where $\Pin(n)^-$ denotes the component which does not contain the identity. We have that $\pm\omega_n$ do not lie in the centre of $\Pin(n)$ (which equals to $\{0,1\}$). In fact, $\omega_n$ commutes with elements in $\text{Spin}(n)$ and
  anti-commutes with elements in $\Pin(n)^-$. This explains that $\pi_1\PO(n)/\pi_0\PO(n)=\{0,1,[\omega_n]\}$ where we  $[\omega_n]$ denotes the class of $\omega_n\in\pi_1\PO(n)$ in $\pi_1\PO(n)/\pi_0\PO(n)$, consisting by $\pm\omega$.
  \begin{enumerate}
  \item Suppose now that $X$ is a closed oriented surface of genus $g$. 
  
  Hence, $\PO(n)$-bundles on $X$ which reduce to $\PSO(n)$ (so with $\mu_1=0$) are classified by 
  $$H^2(X,\pi_1\PO(n))/\Z_2\cong\pi_1\PO(n)/\Z_2\cong\{0,1,[\omega_n]\}.$$

If $\mu_1\neq 0$, then reasoning in a similar manner to the $\Or(2)$ case in \eqref{eq:H2 O2}, one shows that $H^2(X,(\Z_4)_{\mu_1})\cong\{[0],[\omega_n]\}\cong\Z_2$ if $n=2\text{ mod }4$. So $0$ and $1$ are identified in this twisted cohomology.
In addition, $0$ and $\omega_n$ are not in the same orbit under the $\pi_0\PO(n)$-action, so the action of $\Z_2=\pi_0\PO(n)$ in $\Z_2=H^2(X,(\Z_4)_{\mu_1})$, induced by \eqref{action pi0G on C2}, is trivial. So $H^2(X,(\Z_4)_{\mu_1})/\Z_2\cong\Z_2$ for $\mu_1\neq 0$ and $n=2\text{ mod }4$. The same argument gives the same result for $n$ multiple of $4$.
Concluding, if $n\geq 4$ is even, $\PO(n)$-bundles over the surface $X$ are classified by characteristic classes 
$$(\mu_1,\mu_2)\in\left(\{0\}\times\{0,1,[\omega_n]\}\right)\cup\left((\Z_2^{2g}\setminus\{0\})\times\Z_2\right).$$
So we see that, when $n\geq 4$ is even, there are precisely $2^{2g+1}+1$ isomorphism classes of topological $\PO(n)$-bundles over such surface. These invariants have interpretations in terms of obstruction theory of bundles \cite[Proposition 3.2]{oliveira:2011}. Namely, a $\PO(n)$-bundle with $\mu_1=0$ (so actually a $\PSO(n)$-bundle) lifts to $\mathrm{SO}(n)$-bundle if and only if $\mu_2=0,1$ and lifts to $\mathrm{Spin}(n)$-bundle if and only if $\mu_2=0$. A $\PO(n)$-bundle with $\mu_1\neq 0$ lifts to $\Or(n)$-bundle if and only it lifts to a $\mathrm{Pin}(n)$-bundle and this happens if and only if $\mu_2=0$. Indeed, if a $\PO(n)$-bundle with $\mu_1\neq 0$ lifts to $\Or(n)$, then one can always choose such a lift to have vanishing second Stiefel-Whitney class; see Remark 3.6 of \cite{oliveira:2011}.

\item Suppose now that $X$ is closed and non-orientable. Then 
$$H^2(X,\pi_1\PO(n))\cong\pi_1\PO(n)/2\pi_1\PO(n).$$ Here, $2\pi_1\PO(n)=\{0\}$, if $n=0\ \text{mod}\ 4$, while $2\pi_1\PO(n)=\{0,1\}$ for $n=2\ \text{mod}\ 4$, 
  so $$H^2(X,\pi_1\PO(n))\cong\begin{cases}
    \Z_2\times\Z_2 & \text{if } n=0\ \text{mod}\ 4 \\
    \{[0],[\omega_n]\}\cong\Z_2 & \text{if } n=2\ \text{mod}\ 4.
  \end{cases}$$ By acting further by $\Z_2=\pi_0\PO(n)$ as above, we see that $\PO(n)$-bundles on the non-orientable surface $X$ which reduce to $\PSO(n)$, so with $\mu_1=0$, are classified by 
  $$H^2(X,\pi_1\PO(n))/\Z_2\cong\begin{cases}
    \{0,1,[\omega_n]\} & \text{if } n=0\ \text{mod}\ 4 \\
    \{[0],[\omega_n]\}\cong\Z_2 & \text{if } n=2\ \text{mod}\ 4.
  \end{cases}$$
 Consider now $\mu_1\neq 0$ and, as above, let $w_X:\pi_1X\to\Aut(\Z)\cong\Z_2$ be the non-trivial orientation character, with $\Z_{w_X}$ the corresponding local system on $X$.
Then, as in \eqref{eq:H2nonorienmu1neq0} and in Example 5(a),
\begin{equation*}
\begin{split}
H^2(X,\pi_1\PO(n)_{\mu_1})&\cong H_0(X,\Z_{w_X}\otimes_{\Z[\pi_1X]}\pi_1\PO(n)_{\mu_1})\\
&\cong\Z_2\otimes_{\Z}\{[0],[\omega_n]\}\cong\Z_2\otimes_{\Z}\Z_2\cong\Z_2,
\end{split}
\end{equation*}
independently of the residue of $n$ modulo $4$. The next action of $\Z_2=\pi_0\PO(n)$ in $\Z_2\cong H^2(X,\pi_1\PO(n)_{\mu_1})$ is trivial as in Example 5(a).

Concluding, projective orthogonal bundles of rank $n\geq 4$ even, over a non-orientable surface $X$ which is a connected sum of $k$ copies of $\R\mathbb{P}^2$, are topologically classified by
$$(\mu_1,\mu_2)\in\begin{cases}
(\{0\}\times\{0,1,[\omega_n]\})\cup (\Z_2^k\setminus\{0\}\times\Z_2) & \text{if } n=0\ \text{mod}\ 4 \\
\Z_2^k\times\Z_2 & \text{if } n=2\ \text{mod}\ 4. \\
\end{cases}$$ 
\end{enumerate}
%\item Consider now the case of $\PO(2)$-bundles over $X$, which we keep assuming to be oriented. We have $\pi_0\PO(2)\cong\Z_2$, while $\pi_1\PO(2)\cong\Z$, with $\pi_1\Or(2)\cong\Z$ as an index $2$ subgroup.  As in the case of $\Or(2)$, $\pi_0\PO(2)$ acts on $\pi_1\PO(2)$ by changing the sign of the generator. So, as in the orthogonal group case, $\PO(2)$-bundles over $X$ which reduce to $\mathrm{PSO}(2)$-bundles, are classified by $H^2(X,\Z)/\Z_2\cong\Z/\Z_2\cong\Z_{\geq 0}$. 
%For $\PO(2)$-bundles with fixed $\mu_1\neq 0$, then one can see, using the same arguments as in the $\Or(2)$-case, that these are classified by $H^2(X,\Z_{\mu_1})/\Z_2\cong H^2(X,\Z_{\mu_1})/\Z_2\cong\Z_2$. Note however that this characteristic class $\mu_2\in\Z_2$ is not the second Stiefel-Whitney class. Rather, $\mu_2$ measures the obstruction to lift the $\PO(2)$-bundle to $\Or(2)$ and this happens exactly when it lifts to the universal cover, as in the $\PO(n)$-case. 
\end{enumerate}
\end{example}

\begin{remark}
The explicit examples \eqref{ex:O(n)}, \eqref{ex:O(2)} and \eqref{ex:PO(n)} are valid not only for the stated groups but also for any group which is homotopically equivalent to it. For example, respectively, $G=\GL(n,\R)$, $n\geq 3$, $G=\GL(2,\R)$ and $G=\PGL(n,\R)$, with $n\geq 4$ even.
\end{remark}

\section{Reminder of obstruction theory in fibrations}\label{Obstruction theory in fibrations}

Here we briefly sketch the theory of obstructions in fibrations which will be used in the proof of Theorem \ref{thm:class-thm}. Some good references, among several others, on this material are \cite{davis-kirk:2001,steenrod:1999,whitehead:1978}. We always take based spaces.

Consider a Hurewicz (or Serre) fibration $p:E\to Y$ with fibre $F$ and let $f:X\to Y$ be a continuous map.
The group $\pi_1F$ acts on $\pi_kF$ (cf. Section 6.16 of \cite{davis-kirk:2001} or Section 4A of \cite{hatcher:2002}) through free homotopies of based maps $S^k\to F$ and the quotient $\pi_kF/\pi_1F$ is $[S^k,F]$. Assume that $F$ is \emph{simple}, meaning that the action of $\pi_1F$ on $\pi_kF$ is trivial for every $k$ (in the course of the proof of Theorem \ref{thm:class-thm}, we will be in this case). Hence $[S^k,F]=\pi_kF$, for all $k$.

\subsection{Obstruction to extending a partial lift}\label{sec:obs-cochain}

Let $X_1$ be the $1$-skeleton of our two-dimensional $\CW$-complex $X$ and let $$g:X_1\longrightarrow E$$ be a partial lift of $f$:
    $$\begin{tikzcd}
&&F\ar{d}\\
    &&E\ar{d}{p}\\
    X_1\ar[rru, dashed, "g"]\ar[r, hook]&X\ar{r}{f}&Y.
    \end{tikzcd}$$
 
Recall the cochain complex $C^*_{\Z[\pi_1X]}(\widetilde{X},\pi_1F)=\Hom_{\Z[\pi_1X]}(C_*(\widetilde{X}),\pi_1F)$, with $\Z[\pi_1X]$ being the group ring of $\pi_1X$ and $C_*(\widetilde{X})$ the $\Z\pi$-module of chains on the universal cover $\pi:\widetilde X\to X$. The obstruction to extending the partial lift $g:X_1\to E$ to a lift $\tilde{f}:X\to E$ is given by a $2$-cochain, not in $C^*_{\Z[\pi_1X]}(\widetilde{X},\pi_1F)$ but rather in $C^*_{\Z[\pi_1X]}(\widetilde{X},\pi_1\mathcal F_{\alpha_1})$, where $\pi_1\mathcal F_{\alpha_1}$ is a local system where $\pi_1F$ is given the $\Z[\pi_1X]$-module via 
\begin{equation}\label{eq:alpha1}
\alpha_1:\pi_1X\stackrel{f_*}{\longrightarrow}\pi_1Y\longrightarrow\Aut(\pi_1F),
\end{equation}
with $\pi_1Y\to\Aut(\pi_1F)$ induced by the monodromy representation $\pi_1Y\to\Aut^h_*(F)$, where $\Aut^h_*(F)$ denotes the space of based homotopy self-equivalences of $F$; see \cite[Proposition 6.62]{davis-kirk:2001} (here one uses again that $F$ is simple). We will just give an idea of the definition of the mentioned $2$-cochain. For details, see Section 7.10 of \cite{davis-kirk:2001}.

Consider a $2$-cell $\tilde{e}_i^2\in C_*(\widetilde{X})$ of $\widetilde{X}$ and let $h_i:\mathbb D\to\widetilde X$ be the corresponding characteristic map ($\mathbb D\subset\R^2$ denotes the closed unit disc), with $\varphi_i=h_i|_{S^1}:S^1\to\widetilde{X}_1$ the attaching map. 
%$$\begin{tikzcd}
%&&&&F\ar{d}\\
%    &&&&E\ar{d}{p}\\
%    S^1\ar{r}{\varphi_i}&\widetilde{X}_1\ar{r}{\pi|_{\widetilde{X}_1}}&X_1\ar[rru, dashed, "g"]\ar[r, hook]&X\ar{r}{f}&Y.
%    \end{tikzcd}$$
    The composite $p\circ g\circ\pi\circ\varphi_i:S^1\to Y$ is null-homotopic because it equals to $f\circ\pi\circ\varphi_i$ which extends to $\mathbb D\to Y$. By lifting this null-homotopy ($p$ is a fibration) we obtain a homotopy between $g\circ\pi\circ\varphi_i$ and a map $S^1\to F$. So $[g\circ\pi\circ\varphi_i]\in [S^1,F]=\pi_1F$ because $F$ is simple. This element of $\pi_1F$ depends however on the null-homotopy of $f\circ\pi\circ\varphi_i$. To remove this dependence, one has to consider the local system $\pi_1\mathcal{F}_{\alpha_1}$.
    
So consider the $2$-cochain 
\begin{equation}\label{obstruction cochain}
c^2(g)\in C^2_{\Z[\pi_1X]}(\widetilde{X},\pi_1\mathcal F_{\alpha_1})
\end{equation}
whose value on $2$-cells $\tilde{e}_i^2$ is $$c^2(g)(\tilde{e}_i^2)=[g\circ\pi\circ\varphi_i]\in\pi_1\mathcal{F}_{\alpha_1}$$
 and then extend linearly. It is an \emph{obstruction $2$-cochain}, since 
 $g$ extends to a lift $\tilde f:X\to E$ of $f$ if and only if $c^2(g)=0$; cf. Theorem 7.37 (and also Theorem 7.1) of \cite{davis-kirk:2001}. In particular if $F$ is simply-connected, the existence of such lift $\tilde f$ is granted.
%that the composite $S^1\xrightarrow{\varphi_i}\widetilde{X}\xrightarrow{\pi} X\xrightarrow{g}Y$ defines a partial section $s_i:S^1\to g^*Y$, where $g^*Y\to \mathbb D$ is the pullback of the fibration $Y$ to $\mathbb D$ and $g=f\circ\pi\circ h_i$. Since $\mathbb D$ is contractible, $g^*Y$ is fibre homotopically trivial, i.e., there is a fibre-homotopy equivalence between $g^*Y\to \mathbb D$ and $\mathbb D\times F\to \mathbb D$. Therefore $s_i:S^1\to \mathbb D\times F\to F$ (hence $g\circ\pi\circ\varphi_i$) gives rise to an element of $\pi_1F$.

\begin{remark}
Even though we will not need this fact here, the cochain $c^2(g)$ is in fact a cocycle so represents a class in $H^2(X,\pi_1\mathcal{F}_{\alpha_1})$. It turns out that such class vanishes if and only if $g$ can be redefined (relative to the $0$-skeleton) in such a way that it can then be extended to a lift $\tilde f:X\to E$ of $f$.
\end{remark}

% The monodromy action $\pi_1Y\to\Aut^h_*(F)$ induces a bijections $[S^1,F]\to[S^1,F]$ for each element of $\pi_1Z$, but since $\pi_1F$ acts trivially on itself such bijections are actually automorphisms of $\pi_1F$, so we have an action of $\pi_1Y$ on $\pi_1F$ and this determines a local coefficient system over $Y$ with fibre $\pi_1F$. By pulling back it to $X$ via $f$, we obtain a local system over $X$ with fibre $\pi_1F$ corresponding to the action $\alpha_1:\pi_1X\xrightarrow{f_*}\pi_1Y\longrightarrow\Aut(\pi_1F)$ of $\pi_1X$. Denote it by $\pi_1\mathcal{F}_{\alpha_1}$. 

Similarly, if one has $g':X_0\to E$ defined on the $0$-skeleton of $X$ such that $p\circ g=f|_{X_0}$ (obviously, this is always possible), then there is an \emph{obstruction $1$-cochain}
\begin{equation}\label{obstruction 1-cochain}
c^1(g')\in C^1_{{\Z[\pi_1X]}}(\widetilde{X},\pi_0\mathcal{F}_{\alpha_0})
\end{equation}
analogously defined, with $\alpha_0:\pi_1X\stackrel{f_*}{\longrightarrow}\pi_1Y\longrightarrow\Aut(\pi_0F)$ defined analogously to $\alpha_1$. Again, it vanishes precisely when $g'$ extends to a partial lift $g:X_1\to E$ of $f$. In particular, this holds if $F$ is connected.

\subsection{Obstruction to the existence of vertical homotopies}

As usual, write $I=[0,1]$.

\begin{definition}\label{def vertical homotopy}
Let $p:E\to Y$ be a fibration and $f:X\to Y$ a map. Two lifts $f_0,f_1:X\to E$ of $f$ are \emph{vertically homotopic} if there is a \emph{vertical homotopy} between them, i.e., a homotopy $K:X\times I\to E$ between $f_0$ and $f_1$ such that $p\circ K:X\times I\to Y$ is the constant homotopy of $f$.
\end{definition}

In other words, the homotopy $K$ preserves the fibres of $p$. This is a stronger relation than simply being homotopic: two lifts of $f$ may be homotopic but not vertically homotopic.

Let $f_0$ and $f_1$ be two lifts of $f:X\to Y$ and let $K_1$ be a vertical homotopy between them on $X_1$, that is between their restrictions $f_0|_{X_1}$ and $f_1|_{X_1}$. We can ask whether $f_0$ and $f_1$
are vertically homotopic, i.e., if we can extend $K_1$ to a vertical homotopy on $X$.
There is indeed an obstruction to the existence of such a vertical homotopy $K$, called the \emph{difference $2$-cochain}     
\begin{equation}\label{difference cochain}
d^2(f_0,K_1,f_1)\in C^2_{{\Z[\pi_1X]}}(\widetilde{X},\pi_2\mathcal{F}_{\alpha_2}),
\end{equation}
where \begin{equation}\label{alpha2}
\alpha_2:\pi_1X\stackrel{f_*}{\longrightarrow}\pi_1Y\longrightarrow\Aut(\pi_2F),
\end{equation}
is defined as $\alpha_1$ in \eqref{eq:alpha1}. $d^2(f_0,K_1,f_1)$ is nothing but a slight modification (an increase of dimension) of the obstruction $2$-cochain defined in  \eqref{obstruction cochain}, since it is an obstruction cochain to extend to $X\times I$ the map on the $2$-skeleton of $X\times I$ $$T=f_0\cup K_1\cup f_1:X\times\{0\}\cup X_1\times I\cup X\times\{1\}\longrightarrow E$$ lifting the constant homotopy $H_f:X\times I\to Y$, $H_f(x,t)=f(x)$. 

To define $d^2(f_0,K_1,f_1)$ more explicitly, take any $2$-cell $e_i^2$ of $X$, so that $e_i^2\times I$ is a $3$-cell of $X\times I$. Consider a cell $\tilde{e}_i^2\times I$ of $\widetilde{X}\times I$ mapping to $e_i^2\times I$ by the covering map $\pi\times 1:\widetilde{X}\times I\to X\times I$, and take the corresponding attaching $\varphi_i:S^2\to\widetilde{X}\times\{0\}\cup\widetilde{X}_1\times I\cup\widetilde{X}\times\{1\}$ be the attaching map.
% So, 
%  $$\begin{tikzcd}
%&&&&F\ar{d}\\
%    &&&&E\ar{d}{p}\\
%    S^2\ar{r}{\varphi_i}&\widetilde{X}\times\{0\}\cup\widetilde{X}_1\times I\cup\widetilde{X}\times\{1\}\ar{r}{\pi\times 1}&X\times\{0\}\cup X_1\times I\cup X\times\{1\}\ar[rru, dashed, "T"]\ar[r, hook]&X\times I\ar{r}{H_f}&Y.
%    \end{tikzcd}$$
Then
\begin{equation}\label{local description of the difference cochain}
d^2(f_0,K_1,f_1)(e_i^2)=c^3(T)(\tilde{e}_i^2\times I)=[T\circ(\pi\times 1)\circ\varphi_i]\in\pi_2\mathcal{F}_{\alpha_2}
\end{equation}
where $c^3(T)$ is the obstruction cochain analogous to \eqref{obstruction cochain} but in dimension $3$. % As above, the reason why $[T\circ(\pi\times 1)\circ\varphi_i]\in\pi_2F$ comes basically from the fact that $\mathbb D\times I$ is contractible. %Notice that, if the vertical homotopy $K_1$ is constant, we can suppose that $T=f_0\cup f_1$, by collapsing $\partial\tilde{e}_i^2\times I$.
We then have that the partial vertical homotopy $T$ extends to a vertical homotopy $X\times I\to E$ if and only if the $2$-cochain $d^2(f_0,K_1,f_1)$ vanishes.

Note finally that, since $X$ has dimension $2$, $d^2(f_0,K_1,f_1)$ is actually a cocycle, therefore defines cohomology class
\begin{equation}\label{eq:differencecocycle}
\delta^2(f_0,K_1,f_1)\in H^2(X,\pi_2\mathcal{F}_{\alpha_2}).
\end{equation}

Again, there is also a similar \emph{difference $1$-cochain} in dimension $1$, 
\begin{equation}\label{difference cochain on dim 1}
d^1(f'_0,K_0,f'_1)\in C^1_{{\Z[\pi_1X]}}(\widetilde{X},\pi_1\mathcal{F}_{\alpha_1})
\end{equation}
which is the obstruction to the existence of vertical homotopies of lifts $f_0',f_1':X_1\to E$ of $f|_{X_1}$, when we are given a partial homotopy, over the $0$-skeleton $X_0$ of $X$, $K_0:X_0\times I\to E$ between $f_0'|_{X_0}$ and $f_1'|_{X_0}$. So $d^1(f'_0,K_0,f'_1)=0$ if and only if $T'=f'_0\cup K_0\cup f'_1:X\times\{0\}\cup X_0\times I\cup X\times\{1\}\longrightarrow E$ can be extended to a vertical homotopy $X_1\times I\to E$ between $f_0',f_1'$. In particular, this is always possible if $F$ is simply-connected.

\section{A reminder of Postnikov sections}

The material in this section is also well-known. We provide references in due course. 

\subsection{Universal fibrations of Eilenberg-Maclane spaces}\label{sec:4.1}

We work within the category of compactly generated based spaces. Let $A$ be an abelian group and, for $n\geq 1$, let $K(A,n)$ be the corresponding Eilenberg-MacLane space.

There is a specific model for $K(A,n)$ which is a topological abelian group. It is constructed by taking the geometrical realization of the simplicial abelian group associated, via the Dold-Kan correspondence \cite{goerss-jardine:1999}, to the chain complex $A[-n]$ which has $A$ in dimension $n$ and $0$ elsewhere. In all that follows we use this model for $K(A,n)$ and we choose for base point the zero of its group structure.

Consider the space $\Map(X,K(A,n))$ of all continuous maps from $X$ to $K(A,n)$, with the topology determined by the compact-open topology, and so that $\Map(X,K(A,n))$ is compactly generated (see Definition 6.2 and Section 6.1.3 of \cite{davis-kirk:2001}). %Define similarly the space $\Map_*(X,K(A,n))$ of based continuous maps from $X$ to $K(A,n)$, with base point the constant map from $X$ to the base point of $K(A,n)$.
For the given model of $K(A,n)$, $\Map(X,K(A,n))$ is a topological abelian group.

A \emph{group-like monoid} is a topological monoid such that the monoid of its path connected components is a group (the difference from a topological group is the possible non-existence of inverses).

Let $\Aut^h(K(A,n))$ denote the group-like monoid of all self-homotopy equivalences of $K(A,n)$. The model of $K(A,n)$ is functorial, thus a group homomorphism $A\to B$ induces a continuous map $K(A,n)\to K(B,n)$ which is also a group homomorphism. It follows that an automorphism of $A$ induces an automorphism of the topological abelian group $K(A,n)$. Hence $\Aut(A)$ acts on $K(A,n)$ by group automorphisms, so we can take their semidirect product $\Aut(A)\ltimes K(A,n)$. This produces an inclusion
\begin{equation}\label{incl-non-pointed-auto}
\Aut(A)\subset\Aut^h(K(A,n))
\end{equation}
given by the composite $\Aut(A)\hookrightarrow \Aut(A)\ltimes K(A,n)\hookrightarrow \Aut^h(K(A,n))$ where the second map is the inclusion
\begin{equation}\label{nu}
 \nu:\Aut(A)\ltimes K(A,n)\hookrightarrow \Aut^h(K(A,n)),\ \ \ \nu(\varphi,x)(z)=\varphi(z)+x,
\end{equation}
 with $+$ denoting the group operation in $A$.
 It turns out that the inclusion \eqref{nu} is a weak homotopy equivalence.

%Recall that, if $G$ is a topological group,  $BG$ denotes its classifying space and $EG\to BG$ is the universal $G$-bundle, with total space $EG$ contractible. From this it follows that the boundary operator (defined by taking as base point the identity of $G$) $\partial:\pi_iBG\stackrel{\cong}{\longrightarrow}\pi_{i-1}G$ is an isomorphism, for all $i\geq 1$. 

Now, there exists \cite{stasheff:1963,may:1975} a unique (up to homotopy equivalence) $\CW$-complex $B\hspace{-0,07cm}\Aut^h(K(A,n))$ which is a \emph{classifying space} for fibre homotopy equivalence classes of fibrations over $X$, with fibre homotopically equivalent to $K(A,n)$. So there is a corresponding \emph{universal fibration} $\Aut^h(K(A,n))\to E\hspace{-0,07cm}\Aut^h(K(A,n))\to B\hspace{-0,07cm}\Aut^h(K(A,n))$, which has similar properties to  the classifying space $BG$ and to the universal $G$-bundle $EG\to BG$, for a topological group $G$. In particular, the set of equivalence classes of the above mentioned fibrations is in bijection with the set of free (i.e.\ non-based) homotopy classes of maps $X\to B\Aut^h(K(A,n))$, i.e., $[X,B\hspace{-0,07cm}\Aut^h(K(A,n))]$. 

Since \eqref{nu} is a weak homotopy equivalence, and since both $B\hspace{-0,07cm}\Aut^h(K(A,n))$ and $B(\Aut(A)\ltimes K(A,n))$ are $\CW$-complexes, it follows that $B\hspace{-0,07cm}\Aut^h(K(A,n))$ is homotopically equivalent to $B(\Aut(A)\ltimes K(A,n))$. In other words, $B(\Aut(A)\ltimes K(A,n))$ is a model for $B\hspace{-0,07cm}\Aut^h(K(A,n))$.

\subsection{Postnikov sections}

Write $$P_kBG$$ for the \emph{$k^{th}$ Postnikov section of $BG$}. By definition (see for example \cite{hatcher:2002} Chapter 4) these $\CW$-complexes fit in a commutative diagram as below, called the \emph{Postnikov tower} of $BG$ $$\xymatrix{&\vdots\ar[d]\\
    &P_2BG\ar[d]\\
    BG\ar[ruu]\ar[ru]\ar[r]&P_1BG}$$ where the vertical maps are fibrations and such that the maps $BG\to P_kBG$ induce isomorphisms $\pi_iBG\cong\pi_iP_kBG$ if $i\leq k$ and, for $i>k$, $\pi_iP_kBG=0$. For each $k$, the homotopy fibre of the fibrations $P_kBG\to P_{k-1}BG$ is the Eilenberg-Maclane space $K(\pi_{k-1}G,k)$.
    
Again, each Postnikov section $P_kBG$ is unique up to homotopy equivalence, hence so is the Postnikov tower. In particular, $P_1BG= B\pi_0G=K(\pi_0G,1)$.

Regarding the second Postnikov section $P_2BG$, there is \cite{robinson:1972} a model for $P_2BG\to B\pi_0G$ which is a fibre bundle with structure group $\pi_0G$. To briefly see this, consider the group-like monoid $\Aut^h(K(\pi_1G,2))$ of self-homotopy equivalences of $K(\pi_1G,2)$. Consider the universal fibration  $$\Aut^h(K(\pi_1G,2))\longrightarrow E\hspace{-0,07cm}\Aut^h(K(\pi_1G,2))\longrightarrow B\hspace{-0,07cm}\Aut^h(K(\pi_1G,2))$$
as defined in the previous section. Recall that $B(\Aut(\pi_1G)\ltimes K(\pi_1G,2))$ is a model for the classifying space  $B\hspace{-0,07cm}\Aut^h(K(\pi_1G,2))$.
Let
\begin{equation}\label{eq:U}
\mathcal{U}=E\hspace{-0,07cm}\Aut^h(K(\pi_1G,2))\times_{\Aut^h(K(\pi_1G,2))}K(\pi_1G,2)\longrightarrow B(\Aut(\pi_1G)\ltimes K(\pi_1G,2))
\end{equation}
be the $K(\pi_1G,2)$-fibration obtained from the above universal fibration.

%Any $K(\pi_1G,2)$-fibration over any finite $\CW$-complex $Z$ is the pullback of $\mathcal{U}$ under a map $Z\to B\hspace{-0,07cm}\Aut^h(K(\pi_1G,2))$ and this produces a bijection between isomorphism classes of such fibrations and the set $[Z,B\hspace{-0,07cm}\Aut^h(K(\pi_1G,2))]$ of homotopy classes of such maps.
 %, thus the latter space is also a classifying space for $K(\pi_1G,2)$-fibrations.
%$[B\pi_0G,B\Aut^h(K(\pi_1G,2))]=[B\pi_0G,B(\Aut(\pi_1G)\ltimes K(\pi_1G,2))]$.

The splitting $$\Aut(\pi_1G)\to \Aut(\pi_1G)\ltimes K(\pi_1G,2),\ \ \varphi\mapsto(\varphi,0)$$ of the exact sequence $K(\pi_1G,2)\hookrightarrow \Aut(\pi_1G)\ltimes K(\pi_1G,2)\to\Aut(\pi_1G)$ yields a section $$e:B\hspace{-0,07cm}\Aut(\pi_1G)\to B(\Aut(\pi_1G)\ltimes K(\pi_1G,2))$$ of the fibration $ B(\Aut(\pi_1G)\ltimes K(\pi_1G,2))\to B\hspace{-0,07cm}\Aut(\pi_1G)$ induced by the above exact sequence. 

Hence the pullback of $\mathcal{U}$ in \eqref{eq:U} under the composition 
\[B\pi_0G\to B\hspace{-0,07cm}\Aut(\pi_1G)\xrightarrow{\ e\ } B(\Aut(\pi_1G)\ltimes K(\pi_1G,2))\] is the fibre bundle $$E\pi_0G\times_{\pi_0G}K(\pi_1G,2)\to B\pi_0G,$$ with fibre $K(\pi_1G,2)$, structure group $\pi_0G$ and where $E\pi_0G\to B\pi_0G$ is the universal $\pi_0G$-principal bundle. Here $\pi_0G$ acts on the $K(\pi_1G,2)$-factor via the composition 
\begin{equation}\label{pi0action}
\pi_0G\xrightarrow{\Psi(-)_*}\Aut(\pi_1G)\hookrightarrow\Aut^h(K(\pi_1G,2)).
\end{equation}
Thus, we use the functorial construction of the Eilenberg-Maclane space $K(\pi_1G,2)$ explained in the preceding  section.
On the other hand, the $\pi_0G$-action on the $E\pi_0G$-factor  is by the monodromy action of $\pi_0G=\pi_1(B\pi_0G)$ (based on the base point of $B\pi_0G=K(\pi_0G,1)$). So, if $\gamma:S^1\to B\pi_0G$ is a loop representing a class in $\pi_0G$ and $v\in E\pi_0G$, then $\gamma\cdot v=v+\tilde\gamma(1)$ where $\tilde\gamma$ is the lift of $\gamma$ through $v$ and $+$ is the group operation on $E\pi_0G$ (recall that $E\pi_0G$ is an abelian group because $B\pi_0G$ is so).

Concluding, using \eqref{pi0action}, if $[a,b]\in E\pi_0G\times_{\pi_0G}K(\pi_1G,2)$ and $g=\gamma\in\pi_0G=\pi_1(B\pi_0G)$, then
\begin{equation}\label{pi0action diagonal}
[a,b]=[a+\tilde\gamma(1),\Psi_*(-g)(b)].
\end{equation}
This is the model for $P_2BG$ we will use: 
\begin{equation}\label{eq:modelforP2BG}
P_2BG=E\pi_0G\times_{\pi_0G}K(\pi_1G,2)\to B\pi_0G.
\end{equation}

\section{The proof}\label{sec: proof}
Now that we have provided the necessary background, we prove the classification Theorem \ref{thm:class-thm}. This will be done by a sequence of lemmas.

\subsection{Replacement of $BG$ by $P_2BG$}

Note that the map $BG\to P_1BG=B\pi_0G$ is $p_{0,*}$ as defined in \eqref{maptoBpi0}. 
The homotopy fibre of the fibration $P_2BG\to P_1BG$ is $K(\pi_1G,2)$, i.e., we obtain the fibration
\begin{equation}\label{Postnikov fibration}
K(\pi_1G,2)\longrightarrow P_2BG\stackrel{p}{\longrightarrow}B\pi_0G
\end{equation}
and the commutative diagram
$$\xymatrix{&K(\pi_1G,2)\ar[d]\\
    &P_2BG\ar[d]^{p}\\
    BG\ar[ru]\ar[r]^{p_{0,*}}&B\pi_0G.}$$

Since the second homotopy group of any Lie group vanishes (cf. \cite{mimura-toda:2000} Chapter VI, Theorem 4.17), we have $\pi_3BG=0=\pi_3P_2BG$ so actually the map $BG\to P_2BG$ induces isomorphisms $\pi_iBG\cong\pi_iP_2BG$, for $i\leq 3$, i.e.,\ it is a $3$-equivalence.  
Since $\dim X=2<3$, we conclude by  \cite{whitehead:1978}, Chapter IV, Theorem $7.16$, that $BG\to P_2BG$ induces a bijection
\begin{equation}\label{P2BG}
 [X,BG]\simeq[X,P_2BG].
\end{equation}
Hence in order to classify $G$-bundles over $X$, it suffices to consider the set $$[X,P_2BG]$$ instead of $[X,BG]$.

As we have seen above, obtaining the classification of $G$-bundles with fixed invariant $\mu_1\in [X,B\pi_0G]$ is equivalent to describing the fibre over $\mu_1$ of the map $\chi:[X,BG]\to[X,B\pi_0G]$ defined in \eqref{maptoBpi0}. 

Now, consider the map induced by the fibration $p$ in \eqref{Postnikov fibration}
\begin{equation}\label{chi}
\chi':[X,P_2BG]\longrightarrow [X,B\pi_0G], \ \ \ \chi'([h])=[p\circ h].
\end{equation}
From \eqref{P2BG} and the previous diagram, we have the commutative diagram 
$$\xymatrix{&[X,P_2BG]\ar[d]^{\chi'}\\
    [X,BG]\ar[ru]\ar[r]^{\chi}&[X,B\pi_0G]}$$
 so the next lemma follows immediately.
\begin{lemma}\label{BG=P2BG}
There is a bijection between the set of isomorphism classes of continuous principal $G$-bundles over $X$ with invariant $\mu_1\in [X,B\pi_0G]$ and the fibre over $\mu_1$ of $\chi'$.
\end{lemma}

Let us now begin applying the obstruction theory sketched in Section \ref{Obstruction theory in fibrations} to our situation, thus to the fibration \eqref{Postnikov fibration} and to the next diagram $$\xymatrix{&K(\pi_1G,2)\ar[d]\\
    &P_2BG\ar[d]^{p}\\
    X\ar@{-->}[ru]\ar[r]^{f}&B\pi_0G}$$
where $f$ is a representative of a class $\mu_1\in[X,B\pi_0G]$.

We want to analyse the fibre $\chi'^{-1}(\mu_1)$ and a first question is if it is non-empty. Being non-empty means that we can lift $f$ to $P_2BG$. To see this, first notice that, since $K(\pi_1G,2)$ is connected, we can always find a partial lift $f'$ of $f$ over the $1$-skeleton, i.e. a lift of $f|_{X_1}$, because \eqref{obstruction 1-cochain} vanishes. In addition, as $K(\pi_1G,2)$ is simply-connected, then the obstruction $2$-cochain \eqref{obstruction cochain} is again zero and hence we can lift $f$ to a map $X\to P_2BG$. So:

\begin{lemma}\label{lemma:surjective}
The map $\chi'$ is surjective.
\end{lemma}

\subsection{Homotopy classes of lifts and vertical homotopy classes of lifts}\label{sec:5.2}

Our next task is then to describe the fibre of $\chi'$ over the class $\mu_1\in[X,B\pi_0G]$. This fibre is in general just a set, namely the set of homotopy classes in $[X,B\pi_0G]$ which project to a representative of the class $\mu_1$, hence to a map homotopic to $f:X\to B\pi_0G$. 
We will describe \emph{all} such homotopy classes of lifts of $f$ in terms of \emph{vertical} homotopy classes of lifts of $f$. Denote the latter set by 
$$[X,P_2BG]_f.$$

First note that if $g:X\to B\pi_0G$ is another representative of the class $\mu_1$, then by lifting a homotopy between $f$ and $g$ we obtain a bijection $[X,P_2BG]_f\cong [X,P_2BG]_g$. 
So we loose nothing by choosing any representative $f$ of $\mu_1$.

Consider the fundamental group $\pi_1(\Map(X,B\pi_0G),f)$, which we regard as the group of homotopy classes of self-homotopies of $f$. 
It acts on $[X,P_2BG]_f$ also by lifting these self-homotopies, that is, if $f_0\in [X,P_2BG]_f$ is a lift of $f$, and $H\in\pi_1(\Map(X,B\pi_0G),f)$ is a self-homotopy of $f$, then $H$ acts on $f_0$ as $H\cdot f_0=\tilde H(-,1)$, where $\tilde H:X\times I\to P_2BG$ is a lift of $H$ such that $\tilde H(-,0)=f_0$. If another lift $\tilde H'$ of $H$ is taken so that $\tilde H'(-,0)=f_0$, then $\tilde H'(-,1)$ is vertically homotopic to $\tilde H(-,1)$ so the action is well-defined.

Since $\pi_0G$ is abelian, $B\pi_0G=K(\pi_0G,1)$ is a topological abelian group, hence so is $\Map(X,B\pi_0G)$. It follows that $\pi_1(\Map(X,B\pi_0G),f)$ is also independent of $f$. 

So the orbit space $$[X,P_2BG]_f/\pi_1(\Map(X,B\pi_0G),f)$$ is independent on the choice of the representative $f$ of $\mu_1$.

\begin{lemma}\label{fibre - lifts of f}
Given $\mu_1\in [X,B\pi_0G]$, there is a bijection of sets $$\chi'^{-1}(\mu_1)\cong[X,P_2BG]_f/\pi_1(\Map(X,B\pi_0G),f),$$
with $f:X\to B\pi_0G$ representing $\mu_1$.
\end{lemma}
\proof 
By definition, $\chi'^{-1}(\mu_1)\subset [X,P_2BG]$. Consider the natural map 
$$[X,P_2BG]_f\to \chi'^{-1}(\mu_1)$$ 
which sends a vertical homotopy class of a lift $g:X\to P_2BG$ of $f$ to its (not necessarily vertical) homotopy class. This is clearly surjective.%: take any class $g$ in $\chi'^{-1}(\mu_1)$ and consider the class in $[X,P_2BG]_f$ of the maps which are equivalent to $g$ by vertical homotopies. The image of such class is the class of $g$ in $\chi'^{-1}(\mu_1)$.

Take two homotopic lifts $f_0$ and $f_1$ of $f$, hence representing the same point in $\chi'^{-1}(\mu_1)$. Projecting the homotopy through $p:P_2BG\to B\pi_0G$, yields a self-homotopy of $f$, which lifts to the homotopy between $f_0$ and $f_1$ we started with. Hence $f_0$ and $f_1$ represent classes in $[X,P_2BG]_f$ which lie in the same orbit of $\pi_1(\Map(X,B\pi_0G),f)$; note that $f_0$ and $f_1$ are vertically homotopic precisely if the induced self homotopy is homotopically trivial. Conversely, by the definition of the action of $\pi_1(\Map(X,B\pi_0G),f)$  by lifting self-homotopies of $f$, we identify the lifts of $f$ which are homotopic. This shows that the above map induces an injective map
$$[X,P_2BG]_f/\pi_1(\Map(X,B\pi_0G),f)\to \chi'^{-1}(\mu_1).$$
Since $[X,P_2BG]_f\to \chi'^{-1}(\mu_1)$ is surjective, then so is the induced map, and we are done.
\endproof

Now we use the model \eqref{eq:modelforP2BG} for $P_2BG$, so that we have the next diagram
\begin{equation}\label{pullback 2}
\xymatrix{&K(\pi_1G,2)\ar[d]\\
    &E\pi_0G\times_{\pi_0G}K(\pi_1G,2)\ar[d]^{p}\\
    X\ar@{-->}[ru]\ar[r]^{f}&B\pi_0G.}
    \end{equation}
The remaining part of the paper will always refer to the above diagram.

%\begin{equation}\label{pullback 2}
 %\xymatrix{&\widetilde{X}\times_{\pi_1X}K(\pi_1G,2)\ar[d]_{q'}\ar[r]&E\pi_0G\times_{\pi_0G}K(\pi_1G,2)\ar[d]^{p'}\ar[r]&Y\ar[d]\\
%&X\ar[r]^(.4)f&B\pi_0G\ar[r]&B\hspace{-0,07cm}\Aut^h(K(\pi_1G,2)).}
%\end{equation}
%Notice that $\pi_1X$ acts on $K(\pi_1G,2)$ via $\Psi(-)_*\circ \mu_{1*}$ defined in \eqref{Psimu1}.

Let $\pi_0G=\pi_1(B\pi_0G)$ act on the set $[X,E\pi_0G\times_{\pi_0G}K(\pi_1G,2)]_f$ of vertical homotopy classes of lifts of $f$ as in \eqref{pi0action}, so via the inclusion $\Aut(\pi_1G)\subset\Aut(K(\pi_1G,2))$. More precisely, a lift $f_0$ of $f$ representing a class in $[X,E\pi_0G\times_{\pi_0G}K(\pi_1G,2)]_f$ can be written as $f_0(x)=[h_1(x),h_2(x)]\in E\pi_0G\times_{\pi_0G}K(\pi_1G,2)$, where the $h_i$ are defined up to the simultaneous $\pi_0G$-action on both factors as in \eqref{pi0action diagonal}. Then $g\in\pi_0G$ acts on the representative $f_0$ as
\begin{equation}\label{action of pi0(G) on vertical homotopy classes}
 (g\cdot f_0)(x)= [h_1(x),\Psi(g)_*(h_2(x))]
\end{equation}
for $x\in X$ and this induces an action on the homotopy classes of lifts of $f$. Indeed $g\cdot f_0$ is another lift of $f$ since $[h_1(x),\Psi(g)_*(h_2(x))]$ and $[h_1(x),h_2(x)]$ lie in the same fibre of $p$.

The next lemma shows that the fibre that we are studying is in bijection with the quotient set for the action \eqref{action of pi0(G) on vertical homotopy classes}.

\begin{lemma}\label{fibre over mu1}
Let $[X,E\pi_0G\times_{\pi_0G}K(\pi_1G,2)]_f$ denote the set of vertical homotopy classes of lifts of $f$ in diagram \eqref{pullback 2}, together with the $\pi_0G$-action \eqref{action of pi0(G) on vertical homotopy classes}.
Then there is a bijection $$\chi'^{-1}(\mu_1)\simeq [X,E\pi_0G\times_{\pi_0G}K(\pi_1G,2)]_f/\pi_0G.$$
\end{lemma}
\proof
Using our model $P_2BG=E\pi_0G\times_{\pi_0G}K(\pi_1G,2)$, we know from Lemma \ref{fibre - lifts of f} that 
\begin{equation}\label{newmodel}
\chi'^{-1}(\mu_1)\simeq [X,E\pi_0G\times_{\pi_0G}K(\pi_1G,2)]_f/\pi_1(\Map(X,B\pi_0G),f),
\end{equation} where the group $\pi_1(\Map(X,B\pi_0G),f)$ acts by lifting self-homotopies of $f$ in \eqref{pullback 2}. 
%Now, since $X$ is $2$-dimensional, Proposition \ref{reduction of structure group of fibration} asserts that the fibrations $f^*P_2BG$ and $\widetilde{X}\times_{\pi_1X}K(\pi_1G,2)$ over $X$ are isomorphic.
%Therefore, lifts of $f$ in \eqref{pullback 1} correspond bijectively to sections of $q$ in \eqref{pullback 1} or equivalently of $q'$ in \eqref{pullback 2} which, in turn, are in bijection with lifts of $f$ in \eqref{pullback 2}. 
%Hence, $$\chi'^{-1}(\mu_1)\simeq [X,E\pi_0G\times_{\pi_0G}K(\pi_1G,2)]_f/\pi_1(\Map(X,B\pi_0G),f)$$ where the action is again by lifting self-homotopies of $f$, in \eqref{pullback 2}.

Now we compute $\pi_1(\Map(X,B\pi_0G),f)$. As we said before, this is independent of $f$, because $\Map(X,B\pi_0G)$ is a topological abelian group. So we compute $\pi_1\Map(X,B\pi_0G)$ taking for base point the constant map $$\text{const}\in\Map(X,B\pi_0G)$$ equal to the base point of $B\pi_0G=K(\pi_0G,1)$. In this case,
\begin{equation}\label{piMap=pi0G}
\begin{split}
\pi_1(\Map(X,B\pi_0G),\text{const})&=[(X\times S^1)/(X\times 1),B\pi_0G]\\ &=[SX\vee S^1,B\pi_0G]\\ &=[S^1,B\pi_0G]\\&=\pi_0G.
\end{split}
\end{equation}
Here, $SX=(X\times I)/(X\times 0\cup X\times 1)$ is the unreduced suspension of $X$, which is simply connected because $X$ is connected, and $\vee$ is the wedge sum. Notice that in \eqref{piMap=pi0G} is not relevant to distinguish classes of based or unbased maps because $B\pi_0G$ is a simple space (it is $1$-simple since $\pi_0G$ is abelian). Moreover, the third equality holds since a homotopy class of a (based) map $f:SX\vee S^1\to B\pi_0G$ is defined by two homotopy classes of (based) maps $f_1:SX\to B\pi_0G$ and $f_2:S^1\to B\pi_0G$, but $[SX,B\pi_0G]=[SX,K(\pi_0G,1)]=H^1(SX,\pi_0G)=0$ because $\pi_1SX=0$.

Denote by $+$ the group operation on $B\pi_0G=K(\pi_0G,1)$ and let $$g\in\pi_0G\cong\pi_1(B\pi_0G).$$ If we represent $g$ by a loop $$\gamma:S^1\longrightarrow B\pi_0G$$ then the self-homotopy of $f$ associated to $g$ is  $H:X\times S^1\to B\pi_0G$ defined by
\begin{equation}\label{self-homotopy of f}
H(x,t)=f(x)+\gamma(t).
\end{equation}
Here, the operation of summing $f$ is the homeomorphism which identifies the component of the constant map $\text{const}:X\to B\pi_0G$ in $\Map(X,B\pi_0G)$ with that of $f$.
Now we use diagram \eqref{pullback 2}.  The element $g\in\pi_0G$ acts on the class of $f_1$ by lifting the self-homotopy $H$ of $f$ given in \eqref{self-homotopy of f} to $E\pi_0G\times_{\pi_0G}K(\pi_1G,2)$ and taking its final value. As $B\pi_0G$ is an abelian group, then $E\pi_0G$ is also an abelian group and we denote again the operation by $+$. A lift of $H$ is given by 
$$\tilde{H}(x,t)=[h_1(x)+\tilde\gamma(t),h_2(x)]$$
where $\tilde\gamma$ is a lift of $\gamma$. 
So, by \eqref{pi0action diagonal}, $$\tilde{H}(x,1)=[h_1(x)+\tilde\gamma(1),h_2(x)]=[h_1(x),\Psi(g)_*(h_2(x))].$$
Hence the action of $\pi_1(\Map(X,B\pi_0G),f)$ by lifting self-homotopies is compatible with the action \eqref{action of pi0(G) on vertical homotopy classes} of $\pi_0G$ and with the identification $\pi_1(\Map(X,B\pi_0G),f)\cong\pi_1(\Map(X,B\pi_0G),\text{const})\cong\pi_0G$ in \eqref{piMap=pi0G}. This, together with \eqref{newmodel}, completes the proof of the lemma.
\endproof

\subsection{Bijection with twisted cohomology}

Consider the fibration $p$ in \eqref{pullback 2}. The corresponding difference cochain \eqref{difference cochain on dim 1} is zero so there are always vertical homotopies between lifts of $f$ on $X_1$.
Thus, given two lifts $f_0$ and $f_1$ of $f$, we can suppose that $f_0|_{X_1}=f_1|_{X_1}$ and that $K_1$ is the constant homotopy between them. We will hence write $$\delta^2(f_0,f_1)$$ instead of $\delta^2(f_0,K_1,f_1)$ for the cocycle \eqref{eq:differencecocycle} determined by the difference cochain \eqref{difference cochain}.
Moreover, the action $\alpha_2$ defined in \eqref{alpha2} becomes $$\alpha_2:\pi_1X\xrightarrow{f_*}\pi_1B\pi_0G\longrightarrow\Aut(\pi_2(K(\pi_1G,2))),$$ that is $$\alpha_2:\pi_1X\xrightarrow{f_*}\pi_0G\longrightarrow\Aut(\pi_1G).$$ But $p$ is a fibre bundle with $\pi_0G$ as the structure group, with $\pi_0G$ acting on the fibre $K(\pi_1G,2)$ as in \eqref{pi0action}, so via $\Psi(-)_*$. Hence the action $\alpha_2$ coincides with $\Psi(-)_*\circ \mu_{1*}$ in \eqref{Psimu1}.
% because the monodromy action $\pi_1BG\to\Aut(\pi_1G)$ in the fibration $G\to EG\to BG$ is sent by $\partial:\pi_1BG\xrightarrow{\cong}\pi_0G$ to the action $\Psi(-)_*$.
Hence we have that $$\delta^2(f_0,f_1)\in H^2(X,\pi_1\mathcal{G}_{\mu_1})$$
with $\pi_1\mathcal{G}_{\mu_1}$ being the local system determined by \eqref{Psimu1}.

The following result (Corollary 6.16 of Chapter VI of \cite{whitehead:1978}) shows that there is a close relation between $[X,E\pi_0G\times_{\pi_0G}K(\pi_1G,2)]_f$ and $H^2(X,\pi_1\mathcal{G}_{\mu_1})$.

\begin{proposition}\label{bijection between vertical homotopy classes of lifts and twisted cohomology group}
Let $f_0:X\to E\pi_0G\times_{\pi_0G}K(\pi_1G,2)$ be any lift of $f:X\to B\pi_0G$. Then the map $$\phi:[X,E\pi_0G\times_{\pi_0G}K(\pi_1G,2)]_f\longrightarrow H^2(X,\pi_1\mathcal{G}_{\mu_1}),\ \ \ \ \phi(f_1)=\delta^2(f_0,f_1)$$ is a bijection.
\end{proposition}

Now we want to see the $\pi_0G$-quotient of Lemma \ref{fibre over mu1} on the side of twisted cohomology $H^2(X,\pi_1\mathcal{G}_{\mu_1})$, via the correspondence of the previous proposition and via the action of $\pi_0G$ in $H^2(X,\pi_1\mathcal{G}_{\mu_1})$ induced from the action in $C^2_{\Z[\pi_1X]}(\widetilde X,\pi_1G)$ defined in \eqref{action pi0G on C2}.

This is the content of the next and final lemma in the proof of Theorem \ref{thm:class-thm}. 

\begin{lemma}\label{final lemma}
The bijection between $[X,E\pi_0G\times_{\pi_0G}K(\pi_1G,2)]_f$ and $H^2(X,\pi_1\mathcal{G}_{\mu_1})$ is $\pi_0G$-equivariant with respect to \eqref{action of pi0(G) on vertical homotopy classes} and \eqref{action pi0G on C2}, respectively. Hence it induces a bijection on the sets of $\pi_0G$-orbits:
\[[X,E\pi_0G\times_{\pi_0G}K(\pi_1G,2)]_f/\pi_0G\simeq H^2(X,\pi_1\mathcal{G}_{\mu_1})/\pi_0G.\]
\end{lemma}
\proof
Recall that we are always dealing with the situation represented by diagram \eqref{pullback 2}.
For the correspondence in Proposition \ref{bijection between vertical homotopy classes of lifts and twisted cohomology group}, one has a choice of a lift $f_0$ of $f$ and the natural one is to choose the zero map, that is, $f_0(x)=[g_1(x),0]$, with $g_1(x)$ such that $p([g_1(x),0])=f(x)$.

Given another lift $f_1$ of $f$, we have already seen that we can suppose that $f_1|_{X_1}=f_0|_{X_1}=0$. With these choices, the cocycle $\delta^2(f_0,f_1)$ is given by \eqref{difference cochain} and \eqref{local description of the difference cochain}, with $T=f_1$, because $f_0$ is the zero map and $f_1|_{X_1}=f_0|_{X_1}=0$ (so $K_1$ is the constant homotopy equal to zero), so that we can collapse $\partial \tilde{e}_i^2\times I\cup\tilde{e}_i^2\times 0$ to one point:
$$\delta^2(f_1)(e_i^2)=[f_1\circ\pi\circ\varphi_i]\in\pi_2K(\pi_1G,2)\cong\pi_1G.$$
On the other hand, by the definition of the action of $\pi_0G$ in $C^2_{\Z[\pi_1X]}(\widetilde{X},\pi_1G)$ in \eqref{action pi0G on C2}, we have that, for $g\in\pi_0G$, 
\begin{equation}\label{action of pi0G on cohomology classes}
 (g\cdot \delta^2(f_1))(e_i^2)=\Psi(g)_*([f_1\circ\pi\circ\varphi_i])\in\pi_1G\cong\pi_2(K(\pi_1G,2)).
\end{equation}
But by \eqref{pi0action}, if we write $f_1(x)=[h_1(x),h_2(x)]$, we conclude that
\begin{equation}\label{action of pi0G on cohomology classes 2}
 \Psi(g)_*([f_1\circ\pi\circ\varphi_i])=[f'_1\circ\pi\circ\varphi_i]
\end{equation} where $f'_1(x)=[h_1(x),\Psi(g)_*(h_2(x))]$.
 
By \eqref{action of pi0(G) on vertical homotopy classes}, \eqref{action of pi0G on cohomology classes} and \eqref{action of pi0G on cohomology classes 2}, we have that $$g\cdot \delta^2(f_1)=\delta^2(g\cdot f_1).$$ We conclude that the bijection stated in Proposition \ref{bijection between vertical homotopy classes of lifts and twisted cohomology group} is $\pi_0G$-equivariant and this completes the proof.
\endproof

Finally, the proof Theorem \ref{thm:class-thm} follows from Lemmas \ref{BG=P2BG}, \ref{lemma:surjective}, \ref{fibre over mu1} and \ref{final lemma}.

\vspace{1cm}

\noindent
      \textbf{André G. Oliveira} \\
      Centro de Matemática da Universidade do Porto, CMUP\\
      Faculdade de Ciências, Universidade do Porto\\
      Rua do Campo Alegre 687, 4169-007 Porto, Portugal\\ 
      \url{www.fc.up.pt}\\
      email: andre.oliveira@fc.up.pt\\
      \url{https://sites.google.com/view/aoliveira}

\vspace{.2cm}
\noindent
\textit{On leave from:}\\
 Departamento de Matemática, Universidade de Trás-os-Montes e Alto Douro, UTAD \\
Quinta dos Prados, 5000-911 Vila Real, Portugal\\ 
\url{www.utad.pt}\\
email: agoliv@utad.pt

\end{document}